\documentstyle[twoside,psfig,epsf,amsfonts]{amsart}

\begin{document}

\title{MODERN STATISTICS BY KRIGING}

\author{T. SUS{\L}O}

\begin{abstract}
We present statistics (S-statistics) based only
on random variable (not random value) with
a mean squared error of mean estimation as a concept
of error.  
\end{abstract}

\maketitle
 
\thispagestyle{empty}

\section{Origin}
\noindent
\label{sec:1}
{\bf Remark.} Notation  
$$\omega^i_j \rho_{ij}$$ 
is equivalent to 
$$
\sum_{i=1}^n \omega^i_j \rho_{ij}=\omega'\rho \ ,
$$ 
where
$$
\begin{array}{cccccccccc}
\omega
&
=
&
\omega_j^i
&
=
&
\underbrace{
\left[
\begin{array}{c}
\omega_j^1 \\
\vdots  \\
\omega_j^n \\
\end{array}
\right]
}_{n\times 1} \ ,
&
\rho
&
=
&
\rho_{ij}
&
=
&
\underbrace{
\left[
\begin{array}{c}
\rho_{1j} \\
\vdots  \\
\rho_{nj} \\
\end{array}
\right] 
}_{n \times 1} \ .
\end{array}
$$
Notation 
$$
\omega^i_j \rho_{ii} \omega^i_j 
$$
is equivalent to
$$
\sum_{i=1}^n\sum_{l=1}^n\omega^i_j \omega^l_j \rho_{il}
=
\sum_{i=1}^n\omega^i_j \sum_{l=1}^n \rho_{il} \omega^l_j
=
\omega' \Lambda \omega \ ,  
$$
where
$$
\begin{array}{ccccccccc}
\rho_{il}
&
=
&
\Lambda
&
=
&
\Lambda'
&
=
&
\rho_{ii}
&
=
&
{\underbrace{
\left[
\begin{array}{ccc}
\rho_{11} & \ldots & \rho_{1n} \\
\vdots  & \ddots & \vdots \\
\rho_{n1}&  \ldots & \rho_{nn} \\
\end{array}
\right]}_{n \times n}} \ . 
\end{array}
$$
Let us consider (e.g. for $j=n+1$) 
variance of the difference $R_j$ of 
two random variables $V_j$ and $\hat{V}_j$, 
where $E\{V_j\}=E\{\hat{V}_j\}=m$, in terms of 
covariance 
$$
\begin{array}{ccc}
D^2\{V_j-\hat{V}_j\}
&
=
&
Cov\{(V_j-\hat{V}_j)(V_j-\hat{V}_j)\} 
\\

&
=
&
Cov\{V_jV_j\} 
-
Cov\{V_j\hat{V}_j\}
-
Cov\{\hat{V}_jV_j\} 
+
Cov\{\hat{V}_j\hat{V}_j\} 
\\

&
=
&
Cov\{V_jV_j\} 
-
2Cov\{\hat{V}_jV_j\} 
+
Cov\{\hat{V}_j\hat{V}_j\}  
\end{array}
$$
introducing the estimation statistics 
$\hat{V}_j
=
\sum_i\omega^i_j V_i
=
\omega^i_j V_i
$ 
$$
\begin{array}{ccc}
D^2\{R_j\}
&
=
&
Cov\{V_jV_j\}-2Cov\{\hat{V}_jV_j\}+Cov\{\hat{V}_j\hat{V}_j\} 
\\

&
=
&
Var\{V_j\}
-2Cov\{\sum_i\omega_j^i V_i V_j\}
+
Cov\{(\sum_i\omega^i_j V_i)(\sum_i\omega^i_j V_i)\} 
\\

&
=
&
\sigma^2-2\sum_i\omega^i_j Cov\{V_i V_j\}
+
\sum_i\sum_l\omega^i_j \omega^l_j Cov\{V_i V_l\} 
\\

&
=
&
\sigma^2
-
2\sigma^2|\omega^i_j \rho_{ij}|
+
\sigma^2 |\omega^i_j \rho_{ii} \omega^i_j |  
\\

&
=
&
\sigma^2
\pm 2\sigma^2\omega^i_j \rho_{ij}
\mp \sigma^2 \omega^i_j \rho_{ii} \omega^i_j  
\end{array}
$$
in terms of correlation function $\rho$
$$
D^2\{R_j\}
=\sigma^2
+2 \sigma^2 
\omega^i_j     
\rho_{ij} 
-\sigma^2 
\omega^i_j 
\rho_{ii} 
\omega^i_j  
$$
if
$$
\omega^i_j 
\rho_{ij} 
< 0      
$$
and
$$
\omega^i_j 
\rho_{ii} 
\omega^i_j < 0 
$$
or
$$
D^2\{R_j\}=
\sigma^2
-2 \sigma^2 \omega^i_j \rho_{ij}     
+\sigma^2 \omega^i_j \rho_{ii} \omega^i_j  
$$
if 
$$
\omega^i_j 
\rho_{ij}  
\ge 0     
$$
and
$$
\omega^i_j 
\rho_{ii} 
\omega^i_j \ge 0 \ . 
$$
The unbiasedness constraint (I condition) 
$$
E\{R_j\}=E\{V_j-\hat{V}_j\}=E\{V_j\}-E\{\hat{V}_j\}=0
$$
is equivalent to
$$
\sum_i \omega^i_j=1 \ . 
$$
The minimization constraint
$$
\frac{\partial D^2\{R_j\}}{\partial \omega^i_j} = 0 \ ,
$$
where
$$
D^2\{R_j\}=\sigma^2
\pm 2 \sigma^2 \omega^i_j \rho_{ij}  
\mp \sigma^2 \omega^i_j \rho_{ii} \omega^i_j 
\mp 2\sigma^2\underbrace{\left(\sum_i\omega^i_j - 1\right)}_0 \mu_j \ , 
$$
produces $n$ equations in the $n+1$ unknowns: 
kriging weights $\omega_j^i$ and a Lagrange 
parameter $\mu_j$ (II condition) 
$$
\begin{array}{cccccl}
{\underbrace{
\left[
\begin{array}{cccc}
\rho_{11} & \ldots & \rho_{1n} &  1 \\
\vdots  & \ddots & \vdots & \vdots \\
\rho_{n1}&  \ldots & \rho_{nn} &  1 \\
\end{array}
\right]}_{n\times(n+1)}}
& 
\cdot 
&
\underbrace{
\left[
\begin{array}{c}
\omega_j^1 \\
\vdots  \\
\omega_j^n \\
\mu_j \\
\end{array}
\right]
}_{(n+1)\times 1}
&
=
&
\underbrace{
\left[
\begin{array}{c}
\rho_{1j} \\
\vdots  \\
\rho_{nj} \\
\end{array}
\right] 
}_{n \times 1} 
\end{array}
$$
multiplied by $\omega^i_j$ 
$$
\omega^i_j \rho_{ii} \omega^i_j + \mu_j \underbrace{\sum_i \omega^i_j}_1
=
\omega^i_j \rho_{ij}
$$
and substituted into $D^2\{R_j\}$
$$
\begin{array}{ccc}
D^2\{R_j\}
&
=
&
E\{[V_j-\hat{V}_j]^2\}-\underbrace{E^2\{V_j-\hat{V}_j\}}_0 \\

&
=
&
E\{[(V_j-m)-(\hat{V}_j-m)]^2\} \\

&
=
&
E\{[V_j-m]^2\}-2(E\{V_j\hat{V}_j\}-m^2)+E\{[\hat{V}_j-m]^2\} \\

&
=
&
\sigma^2
-2 \sigma^2 |\omega^i_j \rho_{ij}|    
+\sigma^2 |\omega^i_j \rho_{ii} \omega^i_j| \\ 

&
=
&
\sigma^2
\pm 2 \sigma^2 \omega^i_j \rho_{ij}   
\mp \sigma^2 \omega^i_j \rho_{ii} \omega^i_j 
\end{array}
$$
give the minimized variance of the field $V_j$ under estimation 
$$
E\{[(V_j-m)-(\hat{V}_j-m)]^2\}
=
\sigma^2
( 1 \pm (\omega^i_j \rho_{ij} + \mu_j) ) 
$$
and these two conditions produce $n+1$ equations
in the $n+1$ unknowns
$$
\begin{array}{cccccl}
{\underbrace{
\left[
\begin{array}{cccc}
\rho_{11} & \ldots & \rho_{1n} &  1 \\
\vdots  & \ddots & \vdots & \vdots \\
\rho_{n1}&  \ldots & \rho_{nn} &  1 \\
1  & \ldots & 1 & 0 \\
\end{array}
\right]}_{(n+1)\times(n+1)}}
& 
\cdot 
&
\underbrace{
\left[
\begin{array}{c}
\omega_j^1 \\
\vdots  \\
\omega_j^n \\
\mu_j \\
\end{array}
\right]
}_{(n+1)\times 1}
&
=
&
\underbrace{
\left[
\begin{array}{c}
\rho_{1j} \\
\vdots  \\
\rho_{nj} \\
1 \\
\end{array}
\right] 
}_{(n+1) \times 1} \ . 
\end{array}
$$

\section{Asymptotic properties of origin}
\noindent
Since
$$
\lim_{n \rightarrow \infty}  
D^2\left\{\sum_{i=1}^n\omega^i_j V_i\right\} 
=
0
$$
then
$$
\lim_{n \rightarrow \infty} \omega^i_j V_i  
=
\lim_{n \rightarrow \infty} \omega^i_j v_i 
$$
and (since)
$$
\lim_{n \rightarrow \infty} E\left\{\sum_{i=1}^n\omega^i_j V_i\right\} 
=
m
$$
then
$$
\lim_{n \rightarrow \infty} E\left\{\sum_{i=1}^n\omega^i_j V_i\right\} 
=
\lim_{n \rightarrow \infty} \sum_{i=1}^n\omega^i_j v_i 
=
m
$$
the minimized variance of the field $V_j$ under estimation
$$
E\{[(V_j-m)-(\hat{V}_j-m)]^2\}
=
\sigma^2
(1\pm(\omega^i_j \rho_{ij} +\mu_j))
$$
has known asymptotic property
$$
\lim_{n \rightarrow \infty} 
E\{[V_j-\hat{V}_j]^2\}  
=
\lim_{n \rightarrow \infty} 
E\{[V_j-\omega^i_j v_i]^2\} 
=
E\{[V_j-m]^2\} 
=
\sigma^2 \ . 
$$

\section{The field under estimation}
\noindent
Let us consider the field $V_j$ under estimation 
$$
\hat{V}_j=\omega^i_j V_i \ ,  
$$
where for auto-estimation holds
$$
\hat{V}_i=\delta^i_i V_i = V_i \ ,  
$$
with minimized variance of the estimation statistics
$$
\begin{array}{ccc}
E\{[\hat{V}_j-m]^2\}
&
=
&
Cov\{(\omega^i_j V_i)(\omega^i_j V_i)\} 
\\

&
=
&
\sum_i\sum_l\omega^i_j \omega^l_j Cov\{V_i V_l\} 
\\

&
=
&
\sigma^2 |\omega^i_j \rho_{ii} \omega^i_j| 
\\

&
=
&
\mp\sigma^2(\omega^i_j \rho_{ij}-\mu_j) \ ,
\end{array}
$$
where  for auto-estimation holds 
$$
E\{[\hat{V}_i-m]^2\}
=
E\{[V_i-m]^2\}
=\sigma^2 
$$
that means outcoming of input value is unknown for mathematical model,
with minimized variance of the field $V_j$ under estimation 
$$
\begin{array}{ccc}
E\{[(V_j-m)-(\hat{V}_j-m)]^2\}
&
=
&
\sigma^2(1\pm(\omega^i_j \rho_{ij} + \mu_j) )    
\end{array}
$$
where for auto-estimation holds 
$$
E\{[(V_i-m)-(\hat{V}_i-m)]^2\}
=
\underbrace{E\{[V_i-m]^2\}}_{\sigma^2}
-
\underbrace{2(E\{V_i\hat{V}_i\}-m^2)}_{2\sigma^2}
+
\underbrace{E\{[\hat{V}_i-m]^2\}}_{\sigma^2} = 0  
$$
that means variance of the field is equal
to variance of the (auto-)estimation statistics (not that
auto-estimation matches observation).

\section{Generalized least squares solution}
\noindent
For $j \rightarrow \infty$
$$
\underbrace{
\left[
\begin{array}{c}
\rho_{1j} \\
\vdots  \\
\rho_{nj} \\
\end{array}
\right] 
}_{n \times 1}
=
\xi
\underbrace{
\left[
\begin{array}{c}
1 \\
\vdots  \\
1 \\
\end{array}
\right] 
}_{n \times 1}
\qquad \xi \rightarrow 0^- ~(\mbox{or} ~\xi \rightarrow 0^+ )
$$ 
and a disjunction of the minimized variance of the field $V_j$ under
estimation 
$$
E\{[V_j-m]^2\}
-
\underbrace{(E\{V_j\hat{V}_j\}-m^2)}_{\mp\sigma^2\xi}
+
\underbrace{E\{\hat{V}_j[\hat{V}_j-V_j]\}}_{\mp\sigma^2\xi}
\quad \mbox{if} \quad
\rho_{ij} \omega^i_j + \mu_j = \xi+ \mu_j=0   
$$
which fulfills its asymptotic property the kriging system   
$$
\begin{array}{cccccl}
{\underbrace{
\left[
\begin{array}{cccc}
\rho_{11} & \ldots & \rho_{1n} & 1 \\
\vdots  & \ddots & \vdots & \vdots \\
\rho_{n1}&  \ldots & \rho_{nn} & 1 \\
1  & \ldots & 1 & 0 \\
\end{array}
\right]}_{(n+1)\times(n+1)}}
& 
\cdot 
&
\underbrace{
\left[
\begin{array}{c}
\omega^1 \\
\vdots  \\
\omega^n \\
- \xi \\
\end{array}
\right]
}_{(n+1)\times 1}
&
=
&
\underbrace{
\left[
\begin{array}{c}
\xi \\
\vdots  \\
\xi \\
1 \\
\end{array}
\right] 
}_{(n+1) \times 1}
&
\end{array}
$$
equivalent to 
$$
\Lambda\omega-\xi F=\xi F 
$$
and
$$
F'\omega=1
$$
where:
$F=(1,\ldots,1)'$, 
$\omega=(\omega^1,\ldots,\omega^n)'$, 
$\Lambda=\Lambda'
=\rho_{ii};~i=1,\ldots,n$, has the least squares solution
$$
\xi=\frac{1}{2 F'\Lambda^{-1}F} 
$$
and
$$
\omega=2\xi\Lambda^{-1}F
=\frac{\Lambda^{-1}F}{F'\Lambda^{-1}F} 
$$
with a mean squared error of mean estimation
$$
E\{[\hat{V}_j-m]^2\}
=
\mp\sigma^2 2\xi \ . 
$$
 
\section{Approach to simple statistics (conclusion)}
\noindent
For white noise 
$$
\begin{array}{ccc}
E\{[(V_j-m)-(\hat{V}_j-m)]^2\}
&
=
&
E\{[V_j-m]^2\}+E\{[\hat{V}_j-m]^2\} \\

&
=
&
\sigma^2(1-\mu_j) \\ 

&
=
&
\sigma^2\left(1+\frac{1}{n}\right)   
\end{array}
$$
since
$$
\lim_{n \rightarrow \infty} E\{[\hat{V}_j-m]^2\}=0
$$
then
$$
\lim_{n \rightarrow \infty} \hat{V}_j 
=
m 
$$
and
$$
\lim_{n \rightarrow \infty}  
E\{[(V_j-m)-(\hat{V}_j-m)]^2\}
=
\sigma^2 \ .
$$

\vspace*{12pt}
\noindent
{\bf Remark.} 
Precession of arithmetic mean can not be identical to $0$ cause 
a straight line fitted to high-noised data by 
ordinary least squares estimator can not have the slope identical to $0$. 
For this reason the estimator of an unknown constant variance
$$
\frac{1}{n}~~\mbox{sum of deviation squares}
$$
in fact is the lower bound for precession of the minimized
variance of the field under estimation 
$$
E\{[(V_j-m)-(\hat{V}_j-m)]^2\}
=
\sigma^2\left(1+\frac{1}{n}\right)   
$$
to increase `a bit' the lower bound for $n \rightarrow \infty$ 
we can effect on weight and reduce total counts $n$ by $1$  
because $1$ is the closest positive integer number to $0$ 
so it is easy to find the closest weight such that 
$$
\frac{1}{n-1} > \frac{1}{n-0} 
$$ 
then the so-called unbiased variance 
$$
\frac{1}{n-1}~~\mbox{sum of deviation squares}
$$
in fact is the simplest estimator of minimized variance of
the field under estimation.

\end{document}